\documentclass[graybox]{svmult}

\usepackage{type1cm}        
\usepackage{makeidx}         
\usepackage{graphicx}        
\usepackage{multicol}        
\usepackage[bottom]{footmisc}
\usepackage{url}

\usepackage{newtxtext}       %
\usepackage{newtxmath}       
\newcommand{\RL}{{\mathcal R}}
\newcommand{\R}{\mathbb{R}}

\usepackage{color}

\makeindex             

\begin{document}

\title*{Several ways to achieve robustness when solving wave propagation problems}
\author{Niall Bootland and Victorita Dolean and Pierre Jolivet and Fr\'ed\'eric Nataf and St\'ephane Operto and Pierre-Henri Tournier}
\authorrunning{Dolean et al.}
\institute{Niall Bootland \at University of Strathclyde, Dept. of Maths and Stats, \email{niall.bootland@strath.ac.uk}
\and Victorita Dolean \at University of Strathclyde, Dept. of Maths and Stats and University C\^ote d'Azur, CNRS, LJAD \email{work@victoritadolean.com}
\and Pierre Jolivet \at University of Toulouse, CNRS, IRIT, \email{pierre.jolivet@enseeiht.fr}
\and Fr\'ed\'eric Nataf \at Sorbonne Universit\'e, CNRS, LJLL, \email{frederic.nataf@sorbonne-universite.fr}
\and St\'ephane Operto \at University C\^ote d'Azur, CNRS, G\'eoazur, \email{stephane.operto@geoazur.unice.fr}
\and Pierre-Henri Tournier \at Sorbonne Universit\'e, CNRS, LJLL \email{tournier@ann.jussieu.fr}}
%
%
\maketitle

\abstract{Wave propagation problems are notoriously difficult to solve. Time-harmonic problems are especially challenging in mid and high frequency regimes. The main reason is the oscillatory nature of solutions, meaning that the number of degrees of freedom after discretisation increases drastically with the wave number, giving rise to large complex-valued problems to solve.
Additional difficulties occur when the problem is defined in a highly heterogeneous medium, as is often the case in realistic physical applications.
For time-discretised problems of Maxwell type, the main challenge remains the significant kernel in curl-conforming spaces, an issue that impacts on the design of robust preconditioners. This has already been addressed theoretically for a homogeneous medium but not yet in the presence of heterogeneities.
In this review we provide a big-picture view of the main difficulties encountered when solving wave propagation problems, from the first step of their discretisation through to their parallel solution using two-level methods, by showing their limitations on a few realistic examples. We also propose a new preconditioner inspired by the idea of subspace decomposition, but based on spectral coarse spaces, for curl-conforming discretisations of Maxwell's equations in heterogeneous media.}

\section{Motivation and challenges}
\label{sec:motivation}
Why do we need robust solution methods for wave propagation problems? Very often in applications, as for example in seismic inversion
, we need to reconstruct the a priori unknown physical properties of an environment from given measurements. From a mathematical point of view, this means solving inverse problems by applying an optimisation algorithm to a misfit functional between the computation and the data. At each iteration of this algorithm we need to solve a boundary value problem involving the Helmholtz equation
\begin{align}
\label{eq:helm}
- \Delta u - \frac{\omega^2}{c^2} u &= f,
\end{align}
where $c = \sqrt{\rho c_P^2}$, $\rho$ is the density of the medium and $c_P$ is the speed of longitudinal waves. Here, $\omega$ is usually given as being the frequency of a localised source and we wish to reconstruct $c = \frac{1}{n}$ from the measured data (here, $n$ is also called the refraction index).

The Helmholtz equation is also known as the reduced wave equation or time-harmonic wave equation.
Solving this equation is mathematically difficult, especially for high wave number $k = \frac{\omega}{c}$, as the solution is oscillatory and becomes more so with increasing $k$. Note that the notion of a high frequency problem is to be understood relative to the size of the computational domain: meaning how many wavelengths are present in the latter. In geophysics, the typically large size of the domain, and therefore the presence of hundreds of wavelengths, makes the problem difficult.

\subsection{Why the time-harmonic problem in mid and high frequency is hard}

\begin{figure}[t]
\sidecaption
\includegraphics[scale=.35]{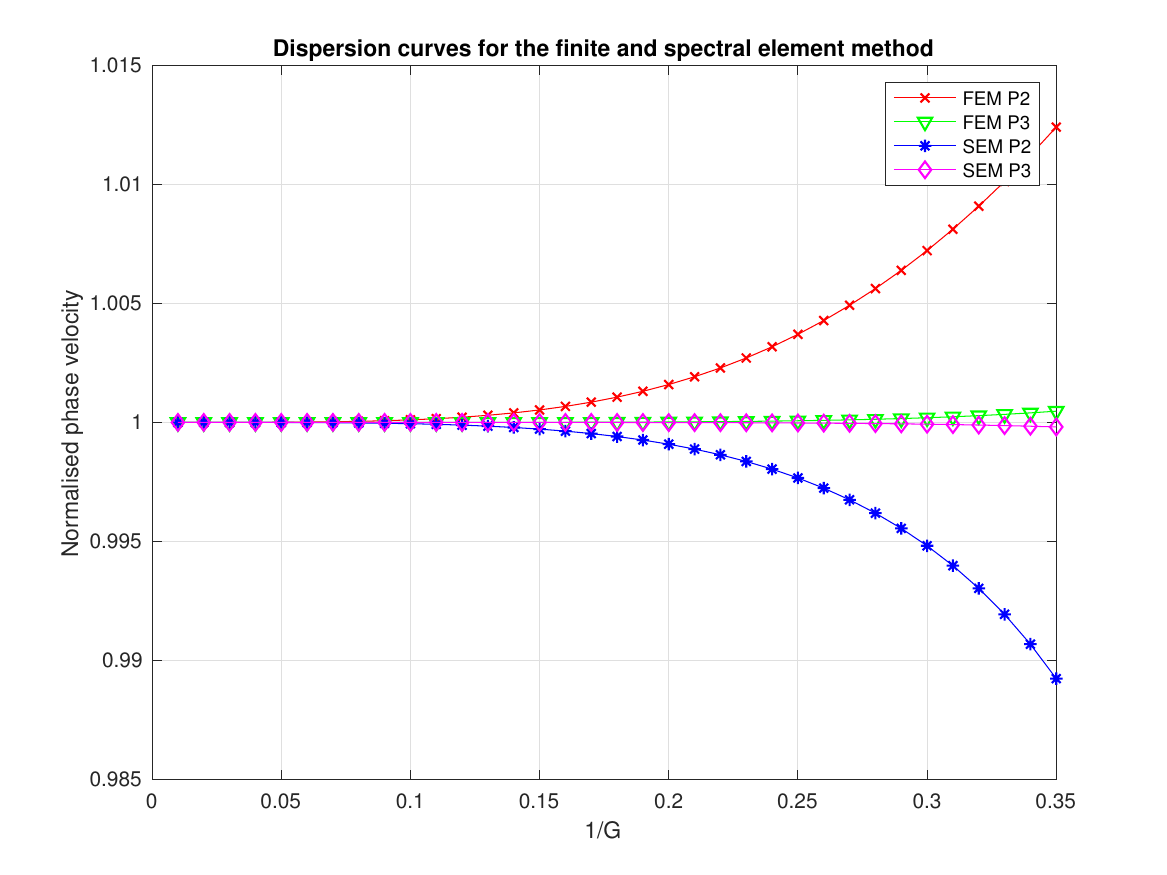}
%
%
\caption{Dispersion curves for finite elements and spectral elements of order 2 and 3: normalised phase velocity (the ratio between the numerical wave speed and the physical one) as a function of the reciprocal number of points per wavelength $\frac{1}{G} = \frac{\omega h}{2\pi}$ for different discretisations. Notice that the use of higher order elements minimises dispersion even for a low value of $G$.}
\label{fig:dispersion}       
\end{figure}

What happens if one wants to approximate this problem with a numerical method? A simple computation in the one-dimensional case shows that the numerical refraction index is different from the physical one and the error depends on the product between the spacing of the grid $h$ and the frequency $\omega$, in other words numerical waves travel at a different speed to physical waves and this is also reflected in the size of error. This is also called the pollution effect and was first highlighted in the seminal paper \cite{Babuska:1997:IPE}. For quasi optimality in the finite element sense we require that $h^p \omega^{p+1}$ be bounded, where $p$ is the order or the precision of the method, as shown in \cite{Melenk:2011:WEC}. To summarise, the high-frequency solution $u$ oscillates at a scale $\frac{1}{\omega}$, therefore the mesh size should be chosen as at least $h \sim \frac{1}{\omega}$ leading to a large number of degrees of freedom. The pollution effect requires $h \ll \frac{1}{\omega}$, namely $h \sim \omega^{-1-\frac{1}{p}}$, therefore in practice one needs an even larger number of degrees of freedom. Note that in order to get a bounded finite element error the constraint is weaker, being {$h \sim \omega^{-1-\frac{1}{2p}}$}, as shown in \cite{Du:2015:PEA}. A trade-off should be found between the number of points per wavelength (ppwl) $G =\frac{\lambda}{h} = \frac{2\pi}{\omega h}$ and the polynomial degree $p$ in order to minimise pollution and this is usually the object of dispersion analysis \cite{Ainsworth:2010:OBS}. This is illustrated in Figure \ref{fig:dispersion}, where we see that the best dispersion properties are achieved when we increase the order of the discretisation or we increase $G$.

Suppose now that we have discretised the equation following the previous rules. We end up with a huge linear system (for a typical application we should expect millions of unknowns) whose size increases with $\omega$ very quickly, especially with more spatial dimensions. The matrix is symmetric and non-Hermitian which makes this system difficult to solve by standard iterative methods, as shown in review paper by Ernst and Gander \cite{Ernst:2012:WDS} or the most recent one by Gander and Zhang \cite{Gander:2019:CIS}. Our aim should be to find the solution in optimal time for large frequencies and our algorithms should not only have good parallel properties but they should also be robust with respect to heterogeneities.

It is well-known that direct solvers, while being robust, have two main drawbacks: their high memory storage and poor parallel properties. On the other hand, iterative methods are not robust but very easy to parallelise. For this reason we consider hybrid methods, such as the naturally parallel compromise of domain decomposition methods, to obtain the best of both worlds. However, how large is truly large? In real applications, problems do not need to be over-resolved (for example, 4 ppwl are enough to perform Full Waveform Inversion with a finite-difference scheme that is specifically tuned to minimise numerical dispersion for this discretisation rule \cite{Amestoy:2016:FFD}) and time-harmonic Helmholtz equations with 50 million degrees of freedom were solved by a parallel direct method \cite{Mary:2017:BLR}. On the other side, when we consider much larger domains (for example via the use of a separate network of nodes rather than with cables) and that the number of nodes is limited, we must switch to iterative or hybrid methods of domain decomposition type. The methods we develop are not only motivated by the current trend in seismic imaging, meaning the development of sparse node devices (OBN) for data acquisition in the oil industry \cite{Blanch:2019:DES}, but in the last decades, since the seminal work of Desprès \cite{Despres:1991:DDM}, they have become the method of choice when solving the discretised Helmholtz equations.

\section{What is the best coarse space for Helmholtz?}
\label{sec:helmholtz}
Consider the decomposition of the computational domain $\Omega$ into $N$ overlapping subdomains $\Omega_j$. The construction of these domains is explained later in Section \ref{seq:cs} and illustrated in Figure \ref{fig:mesh}. We usually solve the system $A\mathbf{u}=\mathbf{b}$ stemming from the finite element discretisation of \eqref{eq:helm} by a preconditioned GMRES method, e.g., in the form $M^{-1}A\mathbf{u}=M^{-1}\mathbf{b}$ with
\begin{align}
M^{-1} &= \sum_{j=1}^N R_j^TD_j B_j^{-1}R_j,
\label{eq:oras}
\end{align}
where $R_j \colon \Omega\rightarrow \Omega_j$ is the restriction operator, $R_j^T \colon \Omega_j\rightarrow \Omega$ the prolongation operator and $D_j$ corresponds to the partition of unity, i.e., it is chosen such that $\sum_{j=1}^N R_j^TD_j R_j = I$. Note also that local matrices $B_j$ are stiffness matrices of local Robin boundary problems
\begin{align*}
(-\Delta - k^2)(u_j) &= f & &\text{in } \Omega_j,\\
\left(\frac{\partial}{\partial n_j}+ ik\right)(u_j) &= 0 & &\text{on } \partial\Omega_j\setminus\partial\Omega.
\end{align*}
We call \eqref{eq:oras} the one-level preconditioner, in particular it is the ORAS preconditioner.

Conventional wisdom in domain decomposition, backed by the definitions of strong and weak scaling, says that one-level preconditioners are not scalable (i.e., their behaviour deteriorates with the number of subdomains $N$). The crucial idea is to add a second level: that is, coarse information that is cheap to compute and immediately available to all subdomains/processors. Suppose that the coarse space is spanned by a matrix Z, then $E= {Z}^*A {Z}$ is the coarse matrix and $H= {Z} E^{-1}{Z}^*$ is the coarse space correction. This coarse space correction can be combined with the one-level preconditioner in an additive or hybrid manner via projectors $P$ and $Q$ ($P = Q = I$ for additive while $P = I - AH$, $Q = I - HA $ provides a hybrid variant)
\begin{align*}
M_2^{-1} = Q M^{-1}P + H.
\end{align*}
This coarse correction can be understood as a solution of a coarser problem on a geometrical grid with a larger spacing for example. For time-harmonic wave propagation problems, the size of the coarse grid is, however, constrained by the wave number. The theory of the grid CS (coarse space) has been introduced by Graham et al.~\cite{Graham:2017:DDP} for a two-level approach to the Helmholtz problem using an equivalent problem with absorption; it has since been extended to the time-harmonic Maxwell equations \cite{Bonazzoli:2019:DDP}. This preconditioner is based on local Dirichlet boundary value problems within the one-level method. An extension to Robin transmission conditions was recently provided in \cite{Graham:2020:DDP}.

The questions we would like to answer are the following: Is the grid coarse space the best choice for heterogeneous problems?
Note also that the definition of the coarse space does not have to be geometrical, we can build more sophisticated coarse spaces based on solving eigenvalue problems. Can we further improve performance by extending the idea of spectral coarse spaces to Helmholtz problems? And, if yes, what kind of modes should be included in the coarse space?

\subsection{Spectral coarse spaces for Helmholtz}
\label{seq:cs}

There are already now a few spectral versions of two-level preconditioners and these are DtN, H-GenEO and $\Delta$-GenEO. For the first two there is no theory available and, while a theory has been developed for the latter, this preconditioner works mainly for low frequency and mildly non-symmetric problems.

The idea of the DtN coarse space was first introduced in \cite{Nataf:2011:ACS} for elliptic problems, further analysed in \cite{Dolean:2012:ATL}, and extended to the Helmholtz equation in \cite{Conen:2014:ACS}. Let $D\subset \Omega$ with internal boundary $\Gamma_D=\partial D\setminus \partial\Omega$ and $v_{\Gamma_D} \colon \Gamma_D\rightarrow \mathbb{C}$. Then the DtN operator is defined as $\mbox{DtN}_D(v_{\Gamma_D}) = \frac{\partial v}{\partial n}|_{\Gamma_D}$
where $v \colon D\rightarrow \mathbb{C}$ is the Helmholtz extension of $v_{\Gamma_D}$ (the solution to a local boundary value problem with Dirichlet value $v_{\Gamma_D}$ on $\Gamma_D$). The DtN coarse space (introduced in \cite{Dolean:2012:ATL}) is based on eigenvalue problems of the DtN operator local to each subdomain:
find $(u_{\Gamma_j},\lambda)\in V(\Gamma_j) \times\mathbb{C}$ such
that
\begin{align*}
\mbox{DtN}_{\Omega_j}(u_{\Gamma_j}) =\lambda u_{\Gamma_j}.
\end{align*}
To provide the {modes} in the coarse space we use the Helmholtz extension $v$. We choose only eigenfunctions with $\lambda$ such that $\mathrm{Re}(\lambda) < k_j$ where \mbox{$k_j = \max_{x\in\Omega_j} k(x)$}. Note that this criterion depends on the local heterogeneity in the problem and is purely heuristic (as explained in detail in \cite{Conen:2014:ACS}). In practice, finding the coarse space vectors amounts to solving local problems depending on Schur complements and mass matrices on the interfaces. By a local Helmholtz extension, we obtain vectors that, after multiplication by the partition of unity and extension by zero, form the matrix $Z$. We do this in each subdomain and combine to give the global coarse space.

The GenEO (Generalised Eigenproblems in the Overlap) coarse space was first developed in \cite{Spillane:2014:ARC} for SPD problems with heterogeneous coefficients, where the heterogeneities do not align with the subdomain decomposition. More precisely, in each $\Omega_j$ we solve {discrete eigenproblems} with local Dirichlet matrices $A_j$ weighted by the partition of unity on one side and the local Neumann matrix $\widetilde{A}_j$ on the other:
\begin{align}
\label{eq:geneo}
D_j A_j D_j u = \lambda \widetilde{A}_j u.
\end{align}
We then choose only eigenfunctions with {eigenvalue $\lambda$} such that $\lambda > \lambda_{\text{min}}$. Note that if we try to replicate this exactly for Helmholtz, the method will fail. For this reason, we need to make some adaptations. The first idea is to use a nearby positive problem to build the coarse space and then use these modes for Helmholtz. This approach is called $\Delta$-GenEO and it is amenable to theory. The second idea is more Helmholtz related in the sense that we only modify the right-hand side of the generalised eigenvalue problem \eqref{eq:geneo} and thus the wave number $k$ is included in the eigenproblem:
\begin{align*}
D_{j} L_{j} D_{j} u = \lambda \widetilde{A}_j u,
\end{align*}
where $L_{j}$ corresponds to the Laplacian part of the problem and $\widetilde{A}_j$ is the Neumann matrix for the Helmholtz operator. Eigenvectors associated to the eigenvalues with $\mathrm{Re}(\lambda) > \lambda_{\text{min}}$ are now those put into the coarse space. We call this method H-GenEO.

\subsection{Comparison of coarse spaces}
\begin{figure}[t]
\centering
\label{fig:marmousi}
\includegraphics[width = \textwidth]{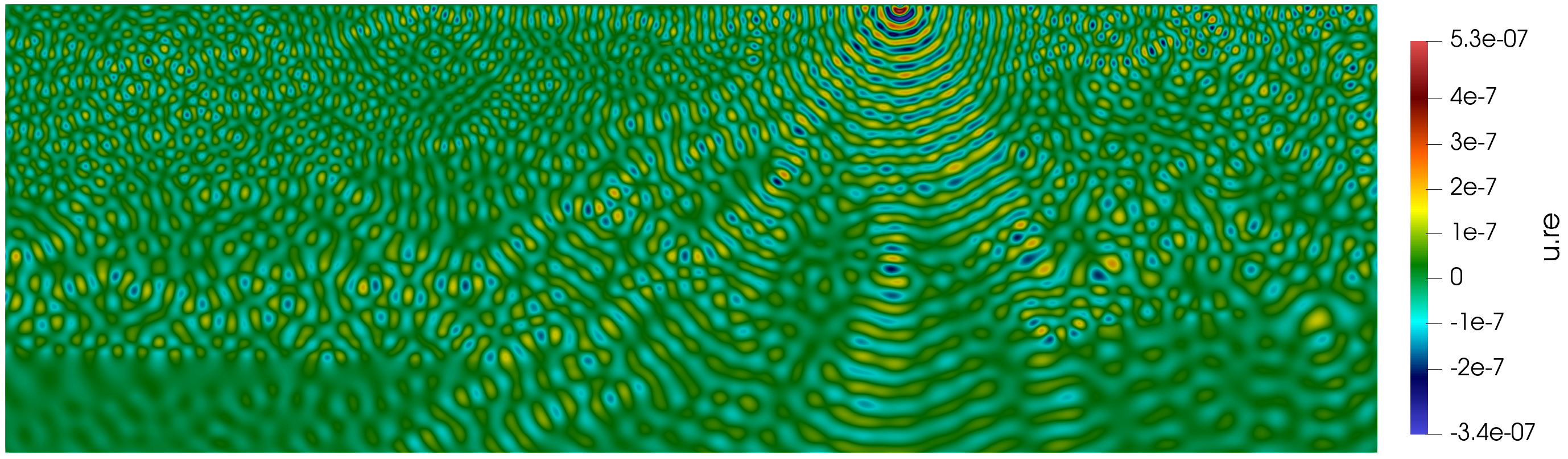}
\caption{The real part of the solution to the Marmousi problem at 20 Hz. The domain is $9.2 \, \mbox{km} \times 3 \, \mbox{km}$.}
\end{figure}

We show how these three two-level methods (grid CS, DtN and H-GenEO) compare on the Marmousi\footnote{ \texttt{https://reproducibility.org/RSF/book/data/marmousi/paper$\_$html/node2.html}} problem \cite{Veersteeg:1994:TME} (see also Figure \ref{fig:marmousi}), which is a 2D geophysical benchmark problem consisting of propagation of seismic waves in a heterogeneous medium from a point source situated towards the surface. This problem is high frequency because of the large number of wavelengths in the domain. For more extensive results and comparative performance tests with these methods on other benchmark problems, see \cite{Bootland:2020:ACS}.

\begin{figure}[t]
\centering
\includegraphics[width = 0.32\textwidth]{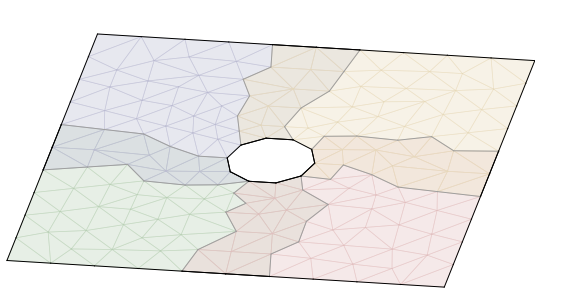}
\includegraphics[width = 0.32\textwidth]{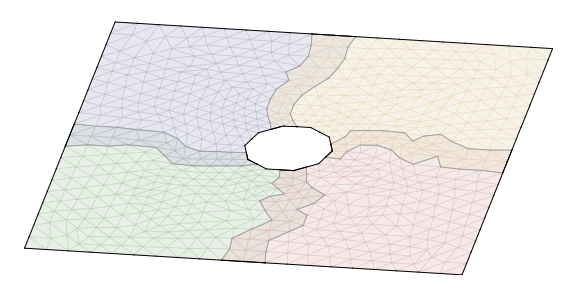}
\includegraphics[width = 0.32\textwidth]{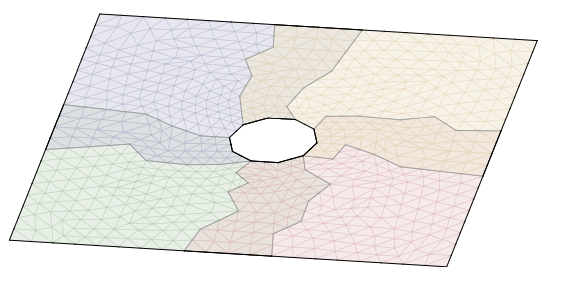}
\caption{Left to right: (a) the coarse mesh, (b) use of minimum overlap (refine the non-overlapping decomposition) and (c) use of coarse overlap (refine the coarse overlapping mesh directly).}
\label{fig:mesh}
\end{figure}

From the practical point of view, a coarse mesh is generated (from which we build the grid coarse space) and this coarse mesh is refined to give the fine mesh; see Figure \ref{fig:mesh}. Alternatively, we can refine on the underlying non-overlapping decomposition and then take minimum overlap. For the discretisation by finite elements (here P2 Lagrange finite elements), we have used FreeFEM. For the domain decomposition and solver we use the FreeFEM library ffddm along with HPDDM and PETSc.\footnote{Software available at \texttt{FreeFem-sources/examples/ffddm} (ffddm) within FreeFEM, \url{https://github.com/hpddm/hpddm} (HPDDM), and \url{https://www.mcs.anl.gov/petsc} (PETSc).}

Note that the grid CS is applied naturally to the configuration (c) (with coarse overlap), whereas for the one-level and spectral methods we can choose between minimum and coarse overlap. In the following tables we report the best results for each method in the most favourable configuration (overlap and number of modes for the spectral coarse spaces). These are the iteration counts of the GMRES method applied to the preconditioned problem with the two-level domain decomposition preconditioner in order to achieve a relative residual tolerance of $10^{-6}$. We consider two scenarios: the {\it under-resolved} case with a discretisation of 5 points per wavelength (Table \ref{Table:Marmousi-P2-nppwl5}) and the {\it over-resolved} case with a discretisation of 10 points per wavelength (Table \ref{Table:Marmousi-P2-nppwl10}) and vary the frequency and the number of subdomains. Low resolution is motivated by applications where high precision is not needed, especially when solving inverse problems by FWI (Full Waveform Inversion). In this case, since the test cases are large, one needs to find a good trade-off between precision and the size of the system to be solved. We refer the reader to the references \cite{Dolean:2020:LSF,Dolean:2021:LSF} where a more extensive numerical study was performed. We notice that in the first scenario the grid CS outperforms the spectral methods (with a slight advantage over the DtN method) whereas in the second scenario the H-GenEO method displays the best performance.

\begin{table}[t]
	\centering
	\caption{Results using the one-level and coarse grid methods for the Marmousi problem when using 5 points per wavelength, varying the frequency $f$ and the number of subdomains $N$.}
	\label{Table:Marmousi-P2-nppwl5}
	\begin{tabular}{cc|ccccc|ccccc|ccccc|ccccc}
	\hline\noalign{\smallskip}
		& & \multicolumn{5}{c|}{One-level} & \multicolumn{5}{c|}{Coarse grid} & \multicolumn{5}{c|}{H-Geneo} & \multicolumn{5}{c}{DtN} \\
		\noalign{\smallskip}\svhline\noalign{\smallskip}
		$f$ & $\text{\#dofs} \ \backslash \ N$ & 10 & 20 & 40 & 80 & 160 & 10 & 20 & 40 & 80 & 160 & 10 & 20 & 40 & 80 & 160 & 10 & 20 & 40 & 80 & 160 \\
		\noalign{\smallskip}\svhline\noalign{\smallskip}
		1 & $4\times10^{3}$ & 26 & 39 & 47 & 64 & $-$ & 15 & 18 & 19 & 20 & $-$ & 9 & 11 & 17 & 21 & $-$ & 6 & 7 & 9 & 6 & $-$ \\
		5 & $1\times10^{5}$ & 53 & 76 & 105 & 154 & 213 & 26 & 29 & 28 & 29 & 31 & 15 & 17 & 26 & 37 & 56 & 7 & 19 & 10 & 8 & 19 \\
		10 & $5\times10^{5}$ & 68 & 102 & 158 & 212 & 302 & 32 & 35 & 41 & 40 & 42 & 33 & 40 & 45 & 56 & 73 & 18 & 19 & 21 & 48 & 29 \\
		20 & $2\times10^{6}$ & 82 & 125 & 178 & 248 & 347 & 34 & 35 & 42 & 43 & 44 & 64 & 83 & 121 & 134 & 157 & 43 & 75 & 77 & 61 & 35 \\
		\noalign{\smallskip}\svhline\noalign{\smallskip}
	\end{tabular}
\end{table}
\begin{table}[t]
	\centering
	\caption{Results using the one-level and coarse grid methods for the Marmousi problem when using 10 points per wavelength, varying the frequency $f$ and the number of subdomains $N$.}
	\label{Table:Marmousi-P2-nppwl10}
	\begin{tabular}{cc|ccccc|ccccc|ccccc|ccccc}
	\hline\noalign{\smallskip}
		& & \multicolumn{5}{c|}{One-level } & \multicolumn{5}{c|}{Coarse grid} & \multicolumn{5}{c|}{H-Geneo } & \multicolumn{5}{c}{DtN} \\
		\noalign{\smallskip}\svhline\noalign{\smallskip}
		$f$ & $\text{\#dofs} \ \backslash \ N$ & 10 & 20 & 40 & 80 & 160 & 10 & 20 & 40 & 80 & 160 & 10 & 20 & 40 & 80 & 160 & 10 & 20 & 40 & 80 & 160\\
		\noalign{\smallskip}\svhline\noalign{\smallskip}
		1 & $2\times10^{4}$ & 30 & 43 & 63 & 97 & $-$ & 16 & 18 & 19 & 21 & $-$ & 7 & 8 & 8 & 13 & $-$ & 4 & 7 & 5 & 6 & $-$\\
		5 & $5\times10^{5}$ & 58 & 87 & 126 & 175 & 246 & 29 & 29 & 34 & 34 & 36 & 10 & 9 & 10 & 10 & 12 & 10 & 11 & 12 & 17 & 24 \\
		10 & $2\times10^{6}$ & 78 & 124 & 172 & 251 & 346 & 35 & 41 & 43 & 46 & 45 & 20 & 16 & 14 & 13 & 13 & 19 & 23 & 25 & 25 & 24 \\
		20 & $8\times10^{6}$ & 92 & 142 & 198 & 272 & 389 & 39 & 47 & 48 & 49 & 49 & 45 & 40 & 34 & 25 & 19 & 35 & 46 & 48 & 56 & 59 \\
		\noalign{\smallskip}\svhline\noalign{\smallskip}
	\end{tabular}
\end{table}
We conclude this comparison by noting that there is no clear advantage in one method over another, all depends on the frequency and precision desired. We have not sought an optimal implementation and the grid CS is the finest possible (which is in principle very expensive), in this sense the timings are not relevant, even if the cost per iteration might be different. In the case of multiple right-hand sides, spectral coarse spaces may have an advantage, although we have not studied this aspect here.

For large-scale geophysical example problems, we have explored extensively the performance of the grid coarse space in \cite{Dolean:2020:IFD,Dolean:2020:LSF}. A few conclusions are stated below:
\begin{itemize}
\item The use of higher order finite elements allow us to minimise dispersion with a minimum number of ppwl, as shown in Figure \ref{fig:dispersion}. A good compromise is the choice of P3 finite elements for which, with 5 ppwl on unstructured meshes, we note a reduction by a factor 2 in the number of degrees of freedom with respect to a finite difference discretisation on uniform meshes.
\item Local solves in domain decomposition methods are usually done by direct methods such as Cholesky factorisation, which is part of the setup phase ahead of the application of the GMRES method. We can already improve performance by replacing the Cholesky method with incomplete Cholesky factorisation.
\item Precision is also important in the parsimony of the computation and the use of single precision highly decreases both the setup and solution times.
\end{itemize}

\section{Can we improve on the auxiliary subspace preconditioner?}
\label{sec:maxwell}

Let us consider the positive (or time-discretised) Maxwell equations
\begin{align*}
\nabla \times (\mu_r^{-1}\nabla\times \mathbf{u} ) + \alpha \varepsilon_r \mathbf{u} &= \mathbf{f} & &\text{in } \Omega,\\
\mathbf{u} \times \mathbf{n} &= 0 & & \text{on } \partial\Omega.
\end{align*}
Here $\mathbf{u}$ is the vector-valued electric field, $\mathbf{f}$ is the source term, $\alpha > 0$ is a constant (e.g., stemming from the time discretisation), and $\mu_r$ and $\varepsilon_r$ are electromagnetic parameters which are uniformly bounded and strictly positive but which we allow to be heterogeneous. We suppose $\Omega$ is a polyhedral computational domain and $\mathbf{n}$ is the outward normal to $\partial\Omega$. After discretisation by N\'ed\'elec elements we obtain
\begin{align}
\label{eq:discretemaxwell}
A\mathbf{U} := (K+\alpha M)\mathbf{U} = \mathbf{b},
\end{align}
where $K\in \mathbb{R}^{n\times n}$ represents the discretisation of the curl--curl operator $\nabla \times (\mu_r^{-1}\nabla\times)$ and $M\in \mathbb{R}^{n\times n}$ is the $\varepsilon_r$-weighted mass matrix computed in the edge element space. Note that matrix $K$ has a huge kernel (all the gradients of $H^1$ functions are part of the kernel of the curl operator) so designing efficient preconditioners for this problem can be challenging.

There is a well-established preconditioner in the literature known as the (nodal) auxiliary space preconditioner (ASP) \cite{Hiptmair:2007:NAS} which is based on a splitting of the space, here $\mathbf{H}(\mathbf{curl})$, by isolating the kernel. The auxiliary space then uses a nodal (Lagrangian) discretisation. The preconditioner is given as
\begin{align*}
{\cal M}^{-1}_{ASP} = \mathrm{diag}( A)^{-1} + P ( \tilde L + \alpha \tilde Q) ^{-1} P^T+\alpha^{-1} C L ^{-1} C^T,
\end{align*}
where $\tilde L$ is the nodal discretisation of the $\mu_r^{-1}$-weighted vector Laplacian operator, $\tilde Q$ is the nodal $\varepsilon_r$-weighted vector mass matrix, $P$ is the matrix form of the nodal interpolation operator between the N\'ed\'elec space and nodal element space, and $C$ is the ``gradient matrix'', which is exactly the null space matrix of $A$ here.

The spectral condition number $\kappa_2({\cal M}^{-1}_{ASP}A)$ of the preconditioned problem is independent of the mesh size but might depend on any heterogeneities present. The natural question is then whether we can improve upon this preconditioner in the case of heterogeneous Maxwell problems?

In order to do this, we use extensively the fictitious space lemma (FSL) of Nepomnyaschickh, which can be considered the Lax--Milgram theorem of domain decomposition \cite{Nepom:1991:DDM}.
\begin{lemma}[Nepomnyaschickh, 1991]
Consider two Hilbert spaces $H$ and $H_D$ along with positive symmetric bilinear forms $a \colon H \times H \rightarrow \R$ and $b \colon H_D \times H_D \rightarrow \R$. The operators $A$ and $B$ are defined as follows
\begin{itemize}
\item $A \colon H \rightarrow H$ such that $(Au,v)=a(u,v)$ for all $u,v\in H$;
\item $B \colon H_D \rightarrow H_D$ such that $(Bu_D,v_D)_D=b(u_D,v_D)$ for all $u_D,v_D\in H_D$.
\end{itemize}
Suppose we have a linear surjective operator $\RL \colon H_D \rightarrow H$ verifying the properties
\begin{itemize}
\item Continuity: $\exists c_R > 0$ such that $\forall u_D\in H_D$ we have
\begin{align*}
a(\RL u_D,\RL u_D) \le c_R b(u_D,u_D).
\end{align*}
\item Stable decomposition: $\exists c_T>0$ such that $\forall u\in H$ $\exists u_D\in H_D$ with $\RL u_D=u$ and
\begin{align*}
c_T b(u_D,u_D) \le a(\RL u_D,\RL u_D)=a(u,u).
\end{align*}
\end{itemize}
Consider the adjoint operator $\RL^* \colon H\rightarrow H_D$ given by
$(\RL u_D, u) = (u_D, \RL^* u)_D$ for all $u_D\in H_D$ and
$u\in H$. Then for all $u\in H$ we have the spectral estimate
\begin{align*}
c_T a(u,u) \le a\left(\RL B^{-1} \RL^* A u,\,u\right) \le c_R a(u,u).
\end{align*}
Thus, the eigenvalues of the preconditioned operator $\RL B^{-1} \RL^* A$ are bounded from below by $c_T$ and from above by $c_R$.
\end{lemma}
In this lemma we have a few ingredients: two Hilbert spaces with the associated scalar products (that are linked by the surjective operator $\RL$) and two symmetric positive bilinear forms. The first of each comes from our problem while the second is for the preconditioner.
Under the assumptions of continuity and stable decomposition, the spectral estimate tells us that the spectral condition number of the preconditioned problem is bounded solely in terms of the constants $c_R$ and $c_T$.

Discretised problems which are perturbations of a singular operator, such as the Maxwell problem in \eqref{eq:discretemaxwell} when $\alpha$ is small, have a huge near-kernel $G \subset \R^{n}$ of $A$, given by the gradient of all $H^1(\Omega)$ functions for example. This near-kernel will be within a space $V_G \subset \R^{n}$, which is the vector space spanned by the sequence $(R_i^T D_i R_i G)_{1\le i \le N}$ so that $G \mathop \subset V_G$. These spaces may not be equal due to the fact that not all the elements of $R_i^T D_i R_i G$ are in $G$, for example, corresponding to the degrees of freedom for which $D_i$ is not locally constant. Nevertheless, since the $D_i$ are related to a partition of unity, we guarantee the inclusion. The space $V_G$ can now serve as a ``free'' coarse space. We denote the coarse space $V_0 := V_G$ and let $Z \in \R^{n_0 \times n}$ be a rectangular matrix whose columns are a basis of $V_0$. The coarse space matrix is then defined in the usual way by $E = Z^TAZ$.

We now need to define all the other ingredients in the FSL. The second Hilbert space is the product space of vectors stemming, for example, from the $n_i$ degrees of freedom on the local subdomains $\Omega_i$ and the $n_0$ coarse space vectors
\begin{align*}
H_D := \R^{n_0} \times \prod_{i=1}^N \R^{n_i}.
\end{align*}
The bilinear form $b$ for the preconditioner is given by the sum of local bilinear forms $b_i$ and the coarse space contribution
\begin{align*}
b(\mathcal{U},\mathcal{V}) &:= (E\mathbf{U}_0,\mathbf{V}_0) + \sum_{i=1}^{N} b_i(\mathbf{U}_i, \mathbf{V}_i), & b_i(\mathbf{U}_i, \mathbf{V}_i) &:= (R_i A R_i^T \mathbf{U}_i, \mathbf{V}_i),
\end{align*}
for $\mathcal{U} = (\mathbf{U}_0, (\mathbf{U}_i)_{1 \le i \le N}) \in H_D$, $\mathcal{V} = (\mathbf{V}_0, (\mathbf{V}_i)_{1 \le i \le N}) \in H_D$.
Finally, the surjective operator $\RL_{AS} \colon H_D \longrightarrow H$ corresponding the additive Schwarz method is given by
\begin{align*}
\RL_{AS}(\mathcal{U}) := Z\mathbf{U}_0 + (I-P_0) \sum_{i=1}^{N} R_i^T \mathbf{U}_i,
\end{align*}
where $P_0$ is the $A$-orthogonal projection on the coarse space $V_0$.
By applying the FSL we obtain a spectral condition number estimate $\kappa (M^{-1}_{AS}A) \le C$, with a bound $C$ that can be large due to heterogeneities in the problem.

How can we improve this preconditioner in this case in order to be robust? We simply build a GenEO space from local generalised eigenproblems in the orthogonal complement of the ``free'' coarse space: find $(\mathbf{V}_{jk},\lambda_{jk}) \in \R^{n_j} \setminus \{0\} \times \R$ such that
\begin{align*}
(I-\xi_{0j}^T) D_j R_j A R_j^T D_j (I-\xi_{0j}) \mathbf{V}_{jk} = \lambda_{jk} \widetilde{A}_j \mathbf{V}_{jk},
\end{align*}
where $\xi_{0j}$ denotes the $b_j$-orthogonal projection from $\R^{n_j}$ on $G_j = R_j G$ and $\widetilde{A}_j$ is the local Neumann matrix for the problem.
We define $V_{j,geneo}^\tau \subset \R^{n}$ to be the vector space spanned by the family of vectors $(R_j^T D_j (I - \xi_{0j}) \mathbf{V}_{jk})_{\lambda_{jk} > \tau}$ corresponding to eigenvalues larger than a chosen threshold parameter $\tau$. Now, collecting over all subdomains $j$, we let $V_{geneo}^\tau$ be the span of all $(V_{j,geneo}^\tau)_{1 \le j \le N}$, which will lead to a new coarse space
\begin{align*}
V_0 := V_G + V_{geneo}^\tau.
\end{align*}
Applying the FSL now yields a spectral condition number estimate of the resulting two-level Schwarz method which is independent of the heterogeneity in the problem.

Several other variants of this approach can be formulated, including with the use of inexact coarse solves in order to more efficiently handle the large coarse space; these theoretical advances can be found in our recent preprint \cite{Bootland:2020:TLD}.

\section{General conclusions}
In this short paper we have offered a brief overview of the main difficulties and some recent solution methods now available to solve Helmholtz equations in the mid and high frequency regimes, which occur in many applications and especially in geophysics. Although there is no established method as the go-to solver, we have proposed a number of different strategies based on two-level domain decomposition methods where the second level comes from the solution of local spectral problems. Indeed, spectral coarse spaces have shown excellent theoretically-proven results for symmetric positive definite problems and currently offer very promising directions to explore for the Helmholtz equation and other wave propagation problems.

The discretisation here is also intertwined with the solution method as solvability and accuracy are very important for wave propagation problems. However, problems in applications do not need to be over-resolved (for example, in full waveform inversion for a discretisation by a finite difference method minimising dispersion 4 ppwl are enough) as this can lead to increasingly large problems whose size is not fully justified by practical reasons. Further, multi-frontal direct solvers based on block low rank approximations have been developed in recent years and problems as large as 50 million unknowns can be tackled successfully by these methods. In this sense, domain decomposition solvers need to be designed with the idea to go beyond these limits while keeping the applicative context in mind.

Last but not least, while not of the same nature, positive Maxwell's equations present different challenges. Here, the auxiliary space preconditioner has successfully been applied to problems where the underlying operator has an infinite dimensional kernel. By exploiting the idea of subspace decomposition together with spectral methods of GenEO type, a new generation of preconditioners, capable of tackling heterogeneous problems, has been introduced. Future work includes an extensive numerical exploration of such an approach on realistic example problems.

Wave propagation problems have been a key source of difficult problems not just for domain decomposition but more widely in scientific computing. As large-scale computing infrastructure continues to evolve and practitioners become ever more ambitious, often driven by industrial challenges, robustness will remain a central theme when designing algorithms for the future. Our work here then contributes some of the most recent ideas towards achieving such desired robustness for domain decomposition methods applied to challenging applications in wave propagation.

\begin{acknowledgement}
The first two authors gratefully acknowledge support from the EPSRC grant EP/S004017/1. The fifth author acknowledges support from the Wind project\footnote{https://www.geoazur.fr/WIND/} funded by Shell, Total and Chevron.
\end{acknowledgement}

%
%

\begin{thebibliography}{99.}%
%
%
\bibitem{Ainsworth:2010:OBS}
M. Ainsworth, H. A. Wajid:  Optimally blended spectral-finite element scheme for wave propagation and nonstandard reduced integration, SIAM J. Numer. Anal., 48(1), pp. 346--371 (2010).

\bibitem{Amestoy:2016:FFD}
P. Amestoy, R. Brossier, A. Buttari, J.-Y. L'Excellent, T. Mary, L. M\'etivier, A. Miniussi, S. Operto: Fast 3D frequency-domain full-waveform inversion with a parallel block low-rank multifrontal direct solver: Application to OBC data from the North Sea, Geophysics, 81(6), pp. R363--R383 (2016).

\bibitem{Babuska:1997:IPE}
I. M. Babu\v{s}ka, S. A. Sauter: Is the pollution effect of the FEM avoidable for the Helmholtz equation considering high wave numbers?, SIAM J. Numer. Anal., 34(6), pp. 2392--2423 (1997).

\bibitem{Blanch:2019:DES}
J. Blanch, J. Jarvis, C. Hurren, Y. Liu, and  L. Hu: Designing an exploration scale {OBN}: {A}cquisition design for subsalt imaging and velocity determination. In {\em {SEG} Technical Program Expanded Abstracts 2019}, pp. 192--196. Society of Exploration Geophysicists (2019).

\bibitem{Bonazzoli:2019:DDP}
M. Bonazzoli, V. Dolean, I. G. Graham, E. A. Spence, P.-H. Tournier: Domain decomposition preconditioning for the high-frequency time-harmonic Maxwell equations with absorption, Math. Comp. 86, pp. 2089--2127 (2017).

\bibitem{Bootland:2020:ACS}
N. Bootland, V. Dolean, P. Jolivet, P.-H. Tournier: A comparison of coarse spaces for Helmholtz problems in the high frequency regime, arXiv preprint \url{arXiv:2012.02678} (2020).

\bibitem{Bootland:2020:TLD}
N. Bootland, V. Dolean, F. Nataf, P.-H. Tournier: Two-level DDM preconditioners for positive Maxwell equations, arXiv preprint \url{arXiv:2012.02388} (2020).

\bibitem{Conen:2014:ACS}
L. Conen, V. Dolean, R. Krause, F. Nataf: A coarse space for heterogeneous Helmholtz problems based on the Dirichlet-to-Neumann operator, J. Comput. Appl. Math., 271, pp. 83--99 (2014).

\bibitem{Despres:1991:DDM}
B. Despr\'es: Domain decomposition method for the Helmholtz problem, C. R. Math. Acad. Sci. Paris. I Math., 311(6), pp. 313--316 (1990).

\bibitem{Dolean:2020:IFD}
V. Dolean, P. Jolivet, P.-H. Tournier, S. Operto: Iterative frequency-domain seismic wave solvers based on multi-level domain-decomposition preconditioners, 82nd EAGE Annual Conference $\&$ Exhibition 2020 (1), pp. 1--5 (2020).

\bibitem{Dolean:2020:LSF}
V. Dolean, P. Jolivet, P.-H. Tournier, S. Operto: Large-scale frequency-domain seismic wave modeling on h-adaptive tetrahedral meshes with iterative solver and multi-level domain-decomposition preconditioners, SEG Technical Program Expanded Abstracts 2020, pp. 2683--2688 (2020).

\bibitem{Dolean:2021:LSF}
V. Dolean, P. Jolivet, P.-H. Tournier, L. Combe, S. Operto, S. Riffo:Large-scale finite- difference and finite-element frequency-domain seismic wave modelling with multi-level domain-decomposition preconditioner, \url{arXiv:2103.14921}, (2021).

\bibitem{Dolean:2012:ATL}
V. Dolean, F. Nataf, R. Scheichl, N Spillane: Analysis of a two-level Schwarz method with coarse spaces based on local Dirichlet--to--Neumann maps, Comput. Methods Appl. Math. 12(4), pp. 391--414 (2012).

\bibitem{Du:2015:PEA}
Y. Du, H. Wu: Preasymptotic error analysis of higher order FEM and CIP-FEM for Helmholtz equation with high wave number, SIAM J. Numer. Anal., 53(2), pp. 782--804 (2015).

\bibitem{Ernst:2012:WDS}
O. G. Ernst., M. J. Gander: Why it is difficult to solve Helmholtz problems with classical iterative methods. In: Graham I., Hou T., Lakkis O., Scheichl R. (eds) Numerical Analysis of Multiscale Problems. LNCSE, vol 83. Springer, Berlin, Heidelberg (2012).

\bibitem{Gander:2019:CIS}
M. J. Gander, H. Zhang: A class of iterative solvers for the Helmholtz equation: Factorizations, sweeping preconditioners, source transfer, single layer potentials, polarized traces, and optimized Schwarz methods, SIAM Rev., 61(1), pp. 3--76 (2019).

\bibitem{Graham:2017:DDP}
I. G. Graham, E. A. Spence, E. Vainikko: Domain decomposition preconditioning for high-frequency Helmholtz problems with absorption, Math. Comp. 86, pp. 2089--2127 (2017).

\bibitem{Graham:2020:DDP}
I. G. Graham, E. A. Spence, J. Zou: Domain Decomposition with local impedance conditions for the Helmholtz equation with absorption, SIAM J. Numer. Anal., 58(5), pp. 2515--2543 (2020).

\bibitem{Hiptmair:2007:NAS}
R. Hiptmair, J. Xu: Nodal auxiliary space preconditioning in {${\bf H}({\bf curl})$} and {${\bf H}({\rm div})$} spaces, SIAM J. Numer. Anal., 45(6), pp. 2483--2509 (2007).

\bibitem{Mary:2017:BLR}
T. Mary: Block low-rank multifrontal solvers: complexity, performance and scalability. PhD Thesis, Universit\'e de Toulouse (2017).

\bibitem{Melenk:2011:WEC}
J. M. Melenk, S. A. Sauter: Wavenumber explicit convergence analysis for Galerkin discretizations of the Helmholtz equation, SIAM J. Numer. Anal., 49(3), pp. 1210--1243 (2011).

\bibitem{Nataf:2011:ACS}
F. Nataf, H. Xiang, V. Dolean, N. Spillane: A coarse space construction based on local Dirichlet-to-Neumann maps, SIAM J. Sci. Comput., 33(4), pp. 1623--1642 (2011).

\bibitem{Nepom:1991:DDM}
S. V. Nepomnyaschikh. Mesh theorems of traces, normalizations of function traces and their inversions. Sov. J. Numer. Anal. Math. Modeling, 6(3), pp. 223--242 (1991).

\bibitem{Spillane:2014:ARC}
N. Spillane, V. Dolean, P. Hauret, F. Nataf, C. Pechstein, R. Scheichl: Abstract robust coarse spaces for systems of PDEs via generalized eigenproblems in the overlaps, Numer. Math. 126(4), pp. 741--770 (2014).

\bibitem{Veersteeg:1994:TME}
R. Versteeg: The Marmousi experience: Velocity model determination on a synthetic complex data set, The Leading Edge 13(9), pp. 927--936, doi:10.1190/1.1437051 (1994).

\end{thebibliography}
%

\end{document}